\input amstex
\documentstyle{amsppt}
%----------------------------------------------------------------
% Title:     A counterexample to Khabibullin's conjecture 
%            for integral inequalities.
% Author:    Ruslan Sharipov
% Comments:  AmSTeX, 10 pages, amsppt style
% MSC-class: 26D10, 26D15, 39B62, 47A63
%----------------------------------------------------------------
%           Replacement for output macro definition
%
\catcode`@=11
\redefine\output@{%
  \def\break{\penalty-\@M}\let\par\endgraf
  \ifodd\pageno\global\hoffset=105pt\else\global\hoffset=8pt\fi  
  \shipout\vbox{%
    \ifplain@
      \let\makeheadline\relax \let\makefootline\relax
    \else
      \iffirstpage@ \global\firstpage@false
        \let\rightheadline\frheadline
        \let\leftheadline\flheadline
      \else
        \ifrunheads@ %\let\makefootline\relax
        \else \let\makeheadline\relax
        \fi
      \fi
    \fi
    \makeheadline \pagebody \makefootline}%
  \advancepageno \ifnum\outputpenalty>-\@MM\else\dosupereject\fi
}
\def\Beta{\mathchar"0\hexnumber@\rmfam 42}
\catcode`\@=\active
%----------------------------------------------------------------
\nopagenumbers
\def\negskp{\hskip -2pt}
\def\Kh{\operatorname{Kh}}
\def\blue#1{#1}

\catcode`#=11\def\diez{#}\catcode`#=6
\catcode`_=11\def\podcherkivanie{_}\catcode`_=8
%\catcode`~=11\def\volna{~}\catcode`~=\active
\def\mycite#1{\cite{\blue{#1}}\immediate\special{ps:
     ShrHPSdict begin /ShrBORDERthickness 0 def}}

\def\mytag#1{%
    \tag#1}
\def\mythetag#1{\thetag{\blue{#1}}\immediate\special{ps:
     ShrHPSdict begin /ShrBORDERthickness 0 def}}
\def\myrefno#1{\no#1}
\def\myhref#1#2{\blue{#2}\immediate\special{ps:
     ShrHPSdict begin /ShrBORDERthickness 0 def}}
\def\myEarXivlink{\myhref{http://arXiv.org}{http:/\negskp/arXiv.org}}

\def\mytheorem#1{\csname proclaim\endcsname{Theorem #1}}
\def\mytheoremwithtitle#1#2{\csname proclaim\endcsname{Theorem #1#2}}
\def\mythetheorem#1{\blue{#1}\immediate\special{ps:
     ShrHPSdict begin /ShrBORDERthickness 0 def}}
\def\mylemma#1{\csname proclaim\endcsname{Lemma #1}}
\def\mylemmawithtitle#1#2{\csname proclaim\endcsname{Lemma #1#2}}
\def\mythelemma#1{\blue{#1}\immediate\special{ps:
     ShrHPSdict begin /ShrBORDERthickness 0 def}}
\def\mycorollary#1{\csname proclaim\endcsname{Corollary #1}}

\def\myconjecture#1{\csname proclaim\endcsname{Conjecture #1}}
\def\myconjecturewithtitle#1#2{\csname proclaim\endcsname{Conjecture #1#2}}
\def\mytheconjecture#1{\blue{#1}\immediate\special{ps:
     ShrHPSdict begin /ShrBORDERthickness 0 def}}

%----------------------------------------------------------------
% Cyrillic fonts definition
%\font\eightcyr=wncyr8
%----------------------------------------------------------------
\pagewidth{360pt}
\pageheight{606pt}
\topmatter
\ \vskip -10pt
\title
A counterexample to Khabibullin's conjecture for integral inequalities.
\endtitle
\rightheadtext{A counterexample to Khabibullin's conjecture \dots}
\author
R.~A.~Sharipov
\endauthor
\address 5 Rabochaya street, 450003 Ufa, Russia\newline
\vphantom{a}\kern 12pt Cell Phone: +7(917)476 93 48
\endaddress
\email \vtop to 30pt{\hsize=280pt\noindent
\myhref{mailto:r-sharipov\@mail.ru}
{r-sharipov\@mail.ru}\newline
\myhref{mailto:R\podcherkivanie Sharipov\@ic.bashedu.ru}
{R\_\hskip 1pt Sharipov\@ic.bashedu.ru}\vss}
\endemail
\urladdr
\vtop to 20pt{\hsize=280pt\noindent
\myhref{http://ruslan-sharipov.ucoz.com/}
{http:/\negskp/ruslan-sharipov.ucoz.com}\newline
\myhref{http://www.freetextbooks.narod.ru/}
{http:/\negskp/www.freetextbooks.narod.ru}\newline
\myhref{http://sovlit2.narod.ru/}
{http:/\negskp/sovlit2.narod.ru}\vss}
\endurladdr
\abstract
    Khabibullin's conjecture for integral inequalities has two numeric
parameters $n$ and $\alpha$ in its statement, $n$ being a positive integer
and $\alpha$ being a positive real number. This conjecture is already 
proved in the case where $n>0$ and $0<\alpha\leqslant 1/2$. However, for 
$\alpha>1/2$ it is not always valid. In this paper a counterexample is 
constructed for $n=2$ and $\alpha=2$. Then Khabibullin's conjecture is 
reformulated in a way suitable for all $\alpha>0$.
\endabstract
\endtopmatter
%\loadbold
%\loadeufb
\TagsOnRight
\document

\head
1. Introduction.
\endhead
\myconjecturewithtitle{1.1}{ (Khabibullin)} Let $\alpha>0$ be a positive number 
and let $q=q(t)$ be a continuous function such that $q(t)\geqslant 0$ for all 
$t>0$. Then the inequality 
$$
\hskip -2em
\int\limits^{\,1}_0\left(\,\,\int\limits^{\,1}_x(1-y)^{n-1}
\,\frac{dy}{y}\right)q(tx)\,dx\leqslant t^{\alpha-1}
\mytag{1.1}
$$
fulfilled for all \kern 2pt $0\leqslant t<+\infty$ implies the inequality
$$
\hskip -2em
\int\limits^{+\infty}_0\!\!q(t)\,\ln\Bigl(1+\frac{1}
{t^{\,2\kern 0.2pt\alpha}}\Bigr)\,dt\leqslant
\pi\,\alpha\prod^{n-1}_{k=1}\Bigl(1+\frac{\alpha}{k}\Bigr).
\mytag{1.2}
$$
\endproclaim
     This conjecture was initially formulated in \mycite{1} and \mycite{2}, 
though in some different form. In \mycite{3} it was reformulated in a form 
very close to the above statement. In \mycite{4} it was proved that the 
conjecture~\mytheconjecture{1.1} is valid for $0<\alpha\leqslant 1/2$. 
Another proof of this result is given in \mycite{5}.\par
      Unfortunately, beyond the scope of\/ $0<\alpha\leqslant 1/2$ the
conjecture~\mytheconjecture{1.1} is not always valid. The main goal of this 
paper is to construct a counterexample to the 
conjecture~\mytheconjecture{1.1} for $n=2$ and $\alpha=2$ and reformulate 
this conjecture in a way suitable for all integer $n>0$ and for all 
$\alpha>0$.
\head
2. The kernel function and the transition function. 
\endhead
      The kernel function $A_n(x)$ and the transition function $\varPhi_n(t)$
associated with the conjecture~\mytheconjecture{1.1} were introduced in 
\mycite{5}. The kernel function is defined by the formula
$$
\hskip -2em
A_n(x)=\int\limits^{\,1}_{\!x}(1-y)^n\,\frac{dy}{y}.
\mytag{2.1}
$$
In terms of the kernel function \mythetag{2.1} the inequality \mythetag{1.1} 
is written as
$$
\hskip -2em
\int\limits^{\,1}_0\!A_{n-1}(x)\,q(t\,x)\,dx\leqslant t^{\alpha-1}\text{, 
\ where \ }t\geqslant 0.
\mytag{2.2}
$$
Omitting the inessential value $t=0$ and changing the variable $x$ for the 
variable $y=t\,x$ in the integral \mythetag{2.2}, we transform it to 
$$
\hskip -2em
\int\limits^{\,t}_0\!A_{n-1}(y/t)\,q(y)\,dy\leqslant t^\alpha
\text{, \ where \ }t>0.
\mytag{2.3}
$$\par
      The concept of the transition function $\varPhi_n(t)$ is more 
complicated. It was introduced in \mycite{5} by means of the following 
formula:
$$
\hskip -2em
\varPhi_n(t)=-\frac{d}{dt}\left(\frac{(-t)^{n+1}}{n!}
\,\frac{d^{\,n+1}\varphi}{dt^{n+1}}\right)\!.
\mytag{2.4}
$$
Here $\varphi=\varphi(t)$ is some smooth function of one variable $t$. 
The basic property of the transition function \mythetag{2.4} is its 
interaction with the kernel function \mythetag{2.1}: 
$$
\hskip -2em
\int\limits^{+\infty}_y\!\!\varPhi_n(t)\,A_n(y/t)\,dt=\varphi(y)
\text{\ \ for \ }n\geqslant 0
\mytag{2.5}
$$
(see Lemma~5.1 in \mycite{5}). In order to apply the formula \mythetag{2.5} 
to Khabibullin's conjecture~\mytheconjecture{1.1} we should specify our 
choice of $\varphi(t)$ in \mythetag{2.4}:
$$
\hskip -2em
\varphi(t)=\ln(1+t^{-2\,\alpha})\text{, \ where \ }\alpha>0.
\mytag{2.6}
$$
The function \mythetag{2.6} depends on $\alpha$. This dependence is 
inherited by the transition function \mythetag{2.4}. Therefore we write 
$\varPhi_n=\varPhi_n(\alpha,t)$. Then \mythetag{2.5} is written as
$$
\hskip -2em
\int\limits^{+\infty}_y\!\!\varPhi_{n-1}(\alpha,t)\,A_{n-1}(y/t)\,dt
=\ln(1+y^{-2\,\alpha})\text{\ \ for \ }n\geqslant 1.
\mytag{2.7}
$$
Note that the function \mythetag{2.6} enters the integrand in the left 
hand side of the inequality \mythetag{1.2}. For this reason the the 
formula \mythetag{2.7} is a bridge binding two inequalities \mythetag{1.1} 
and \mythetag{1.2} in Khabibullin's conjecture.
\head
3. Application to Khabibullin's conjecture.
\endhead
     Let's recall that the first inequality of Khabibullin's 
conjecture~\mytheconjecture{1.1} is written as \mythetag{2.3}. Let's
multiply both sides of \mythetag{2.3} by $\varPhi_{n-1}(\alpha,t)$:
$$
\hskip -2em
\int\limits^{\,t}_0\!\varPhi_{n-1}(\alpha,t)
\,A_{n-1}(y/t)\,q(y)\,dy\leqslant \varPhi_{n-1}(\alpha,t)\ t^\alpha.
\mytag{3.1}
$$
The inequality \mythetag{3.1} would follow from \mythetag{2.3} provided
$\varPhi_{n-1}(\alpha,t)\geqslant 0$. But actually, positive values of 
$\varPhi_{n-1}(\alpha,t)$ alternate with negative ones. For this reason
we subdivide the half-line $t>0$ into two subsets $M_{\sssize +}$ and
$M_{\sssize -}$:
$$
\hskip -2em
\aligned
&M_{\sssize +}=M_{\sssize +}(n,\alpha)=\{t\in\Bbb R\!:\ t>0
\text{\ \ and \ }\varPhi_{n-1}(\alpha,t)\geqslant 0\},\\
&M_{\sssize -}=M_{\sssize -}(n,\alpha)=\{t\in\Bbb R\!:\ t>0
\text{\ \ and \ }\varPhi_{n-1}(\alpha,t)<0\}.
\endaligned
\mytag{3.2}
$$
The inequality \mythetag{3.1} is valid for $t\in M_{\sssize +}(n,\alpha)$. 
If\/ $t\in M_{\sssize -}(n,\alpha)$, the inequality \mythetag{3.1} is not 
valid. In this case we recall that the kernel function $A_n(x)$ is positive, 
i\.\,e\. $A_n(x)>0$ for $0<x<1$ (see \mycite{5}). The function $q(y)$ is 
non-negative according to the conjecture~\mytheconjecture{1.1}. As for the 
function $\varPhi_{n-1}(\alpha,t)$, it is negative for 
$t\in M_{\sssize -}(n,\alpha)$ (see \mythetag{3.2}). As a result we get 
the following inequality for $t\in M_{\sssize -}(n,\alpha)$:
$$
\hskip -2em
\int\limits^{\,t}_0\!\varPhi_{n-1}(\alpha,t)
\,A_{n-1}(y/t)\,q(y)\,dy\leqslant 0.
\mytag{3.3}
$$
Combining the inequalities \mythetag{3.1} and \mythetag{3.3}, we write
$$
\int\limits^{\,t}_0\!\varPhi_{n-1}(\alpha,t)
\,A_{n-1}(y/t)\,q(y)\,dy\leqslant
\cases \varPhi_{n-1}(\alpha,t)\ t^\alpha &\text{for \ }t\in M_{\sssize +},
\\\qquad 0 &\text{for \ }t\in M_{\sssize -}.
\endcases
\mytag{3.4}
$$
Now let's integrate \mythetag{3.4} over $t$ from $0$ to infinity:
$$
\hskip -2em
\int\limits^{+\infty}_0\!\left(\,\int\limits^{\,t}_0
\!\varPhi_{n-1}(\alpha,t)\,A_{n-1}(y/t)\,q(y)\,dy\!\right)\!dt
\,\,\leqslant\!\int\limits_{t\in M_{\sssize +}}
\!\!\!\varPhi_{n-1}(\alpha,t)\ t^\alpha\,dt.
\mytag{3.5}
$$
Upon changing the order of integration in the left hand side of 
\mythetag{3.5} we get
$$
\int\limits^{+\infty}_0\left(\int\limits^{+\infty}_y
\!\varPhi_{n-1}(\alpha,t)\,A_{n-1}(y/t)\,dt\!\right)\!q(y)\,dy
\,\,\leqslant\!\int\limits_{t\in M_{\sssize +}}
\!\!\!\varPhi_{n-1}(\alpha,t)\ t^\alpha\,dt.
\quad
\mytag{3.6}
$$
Applying \mythetag{2.7} to \mythetag{3.6}, we derive the following 
inequality:
$$
\hskip -2em
\int\limits^{+\infty}_0\!\ln(1+y^{-2\,\alpha})\,q(y)\,dy
\,\,\leqslant\!\int\limits_{t\in M_{\sssize +}}
\!\!\!\varPhi_{n-1}(\alpha,t)\ t^\alpha\,dt.
\mytag{3.7}
$$\par
     The integral in the left hand side of \mythetag{3.7} coincides 
with that of the second inequality \mythetag{1.2} in Khabibullin's 
conjecture~\mytheconjecture{1.1}. There is the formula 
$$
\hskip -2em
\varPhi_n(\alpha,t)=\frac{4\,\alpha^2}{t}\cdot
\frac{t^{2\,\alpha}\,P_n(\alpha,z)}{(1+t^{2\,\alpha})^{n+2}}.
\mytag{3.8}
$$
It was derived in \mycite{5}. \pagebreak Here $z=t^{2\,\alpha}$ and 
$P_n(\alpha,z)$ is a polynomial of the degree $n$ with respect to 
the variable $z$. From the formula \mythetag{3.8} we derive
$$
\hskip -2em
\varPhi_n(\alpha,t)=
\cases 
O(t^{2\,\alpha-1})&\text{as \ }t\longrightarrow +0,\\
O(t^{-2\,\alpha-1})&\text{as \ }t\longrightarrow +\infty.
\endcases 
\mytag{3.9}
$$
Due to \mythetag{3.9} the integral in the right hand side of 
\mythetag{3.7} is finite. In's value is a positive number depending 
on $n$ and $\alpha$. Let's denote it through $C(n,\alpha)$: 
$$
\hskip -2em
C(n,\alpha)\ =\!\int\limits_{t\in M_{\sssize +}}
\!\!\!\varPhi_{n-1}(\alpha,t)\ t^\alpha\,dt<\infty.
\mytag{3.10}
$$
Apart from \mythetag{3.10}, we consider two other integrals
$$
\xalignat 2
&\hskip -2em
\int\limits_{t\in M_{\sssize -}}
\!\!\!\varPhi_{n-1}(\alpha,t)\ t^\alpha\,dt,
&&\int\limits^{+\infty}_0\!\!\varPhi_{n-1}(\alpha,t)\ t^\alpha\,dt.
\mytag{3.11}
\endxalignat
$$
Due to \mythetag{3.9} both integrals \mythetag{3.11} are finite. 
According to \mythetag{3.2}, the first of them is non-positive. 
These integrals are related to \mythetag{3.10} as follows:
$$
\hskip -2em
C(n,\alpha)\ +\!\!\int\limits_{t\in M_{\sssize -}}
\!\!\!\varPhi_{n-1}(\alpha,t)\ t^\alpha\,dt\ =\int\limits^{+\infty}_0
\!\!\varPhi_{n-1}(\alpha,t)\ t^\alpha\,dt.
\mytag{3.12}
$$
In \mycite{5} the second integral \mythetag{3.11} was calculated 
explicitly. As appears, this integral coincides with the number in 
the right hand side of the second inequality \mythetag{1.2} in 
Khabibullin's conjecture~\mytheconjecture{1.1}, i\.\,e\. we have  
$$
\hskip -2em
\int\limits^{+\infty}_0\!\!\varPhi_{n-1}(\alpha,t)\ t^\alpha\,dt
=\pi\,\alpha\prod^{n-1}_{k=1}\Bigl(1+\frac{\alpha}{k}\Bigr).
\mytag{3.13}
$$
Since the first integral \mythetag{3.11} is non-positive, from 
\mythetag{3.12} and \mythetag{3.13} we derive 
$$
\hskip -2em
0<\pi\,\alpha\prod^{n-1}_{k=1}\Bigl(1+\frac{\alpha}{k}\Bigr)
\leqslant C(n,\alpha)<\infty
\mytag{3.14}
$$
The formula \mythetag{3.10} for $C(n,\alpha)$ is not so simple as 
compared to \mythetag{3.13}. Nevertheless, $C(n,\alpha)$ should be 
considered as a quite certain quantity that can be effectively 
computed numerically for each particular numeric value of $n$ and 
$\alpha$. Replacing the right hand side of the inequality \mythetag{1.2} 
by $C(n,\alpha)$, we cam formulate Khabibullin's 
conjecture~\mytheconjecture{1.1} not as a conjecture, but as a proved 
result.
\mytheorem{3.1} Let $\alpha>0$ be a positive number and let $q=q(t)$ 
be a continuous function such that $q(t)\geqslant 0$ for all $t>0$. 
Then the inequality 
$$
\pagebreak
\int\limits^{\,1}_0\left(\,\,\int\limits^{\,1}_x(1-y)^{n-1}
\,\frac{dy}{y}\right)q(tx)\,dx\leqslant t^{\alpha-1}
$$
fulfilled for all \kern 2pt $0<t<+\infty$ implies the inequality
$$
\hskip -2em
\int\limits^{+\infty}_0 q(t)\,\ln\Bigl(1+\frac{1}
{t^{\,2\kern 0.2pt\alpha}}\Bigr)\,dt\leqslant C(n,\alpha).
\mytag{3.15}
$$
\endproclaim
\head
4. A counterexample to Khabibullin's conjecture.
\endhead
      In this section we show that Khabibullin's conjecture in its 
form \mytheconjecture{1.1} is not valid for $\alpha>1/2$. Since the 
first inequality \mythetag{1.1} in Khabibullin's conjecture for $t>0$ 
is already transformed to \mythetag{2.3}, we introduce the function
$$
\hskip -2em
g(t)=\int\limits^{\,t}_0\!A_{n-1}(y/t)\,q(y)\,dy.
\mytag{4.1}
$$
The formula \mythetag{4.1} is called the direct conversion formula. 
It was studied in \mycite{6}, where the inverse conversion formula 
expressing $q(t)$ back through $g(t)$ was derived:
$$
\hskip -2em
q(t)=\frac{d^{\kern 0.3pt n}}{dt^{n}}\!\left(\frac{t^n
\,g'(t)}{(n-1)!}\right)\!.
\mytag{4.2}
$$
As it was shown in \mycite{6}, under the assumptions of Khabibullin's 
conjecture~\mytheconjecture{1.1} the function \mythetag{4.1} is an $(n+1)$ 
times differentiable function satisfying the inequality
$$
\hskip -2em
0\leqslant g(t)\leqslant t^\alpha\text{\ \ for all \ }t>0.
\mytag{4.3}
$$\par
     Before constructing a counterexample to the 
conjecture~\mytheconjecture{1.1} we specify $n=2$ and $\alpha=2$. Then 
the transition function \mythetag{3.8} turns to 
$$
\hskip -2em
\varPhi_{n-1}(\alpha,t)=\frac{16\,t^3\,P_1(\alpha,z)}{(1+t^4)^3},
\mytag{4.4}
$$
where $z=t^4$ and $P_1(\alpha,z)$ is the following polynomial:
$$
\hskip -2em
P_1(\alpha,z)=(2\,\alpha+1)\,z+(1-2\,\alpha)=5\,z-3
\mytag{4.5}
$$
(see the formula \thetag{7.6} in \mycite{5}). Substituting \mythetag{4.5} 
into \mythetag{4.4}, we derive 
$$
\hskip -2em
\varPhi_1(2,t)=\frac{16\,t^3\,(5\,t^4-3)}{(1+t^4)^3}.
\mytag{4.6}
$$
The basic feature of the function \mythetag{4.6} is that it is not 
always positive. According to \mythetag{3.2}, we have two non-empty 
subsets $M_{\sssize -}$ and $M_{\sssize +}$ of the half-line $t>0$:
$$
\xalignat 2
&\hskip -2em
M_{\sssize -}=\{t\in\Bbb R\!:\ 0<t<t_0\},
&&M_{\sssize +}=\{t\in\Bbb R\!:\ t\geqslant t_0\},
\quad
\mytag{4.7}
\endxalignat
$$
where $t_0=\root 4\of{3/5}$ is the root of the polynomial $5\,t^4-3$ 
in the numerator of the fraction in \mythetag{4.6}. \pagebreak Relying 
on \mythetag{4.7}, we construct the function $g=g(t)$ defining it by 
two different formulas in $M_{\sssize -}$ and in $M_{\sssize +}$. 
For this purpose we define the polynomial
$$
\hskip -2em
h(t)=\frac{(t-t_0)^4}{{t_0}^{\kern -3pt 4}}.
\mytag{4.8}
$$
Since $\alpha=2$, we set $g(t)=t^\alpha=t^2$ for $t\in M_{\sssize +}$. 
For $t\in M_{\sssize -}$ we define $g(t)$ by a spline polynomial 
composed with the use of the polynomial \mythetag{4.8}: 
$$
\hskip -2em
g(t)=\cases t^2\,(1-\varepsilon h(t))&\text{for \ }0<t<t_0,\\
\kern 2.5em t^2&\text{for \ }t\geqslant t_0.\\
\endcases
\mytag{4.9}
$$
Let's recall that $n=2$ and $n+1=3$. Hence \mythetag{4.9} should be a 
three times differen\-tiable function. This leads to the following 
conditions at the point $t=t_0$:
$$
\xalignat 2
&\hskip -2em
\lim_{t\to t_0}g(t)={t_0}^{\kern -3pt 2},
&&\lim_{t\to t_0}g'(t)=2\,t_0,\\
\vspace{-1ex}
&&&\mytag{4.10}\\
\vspace{-1ex}
&\hskip -2em
\lim_{t\to t_0}g''(t)=2,
&&\lim_{t\to t_0}g'''(t)=0.
\endxalignat
$$
It is easy to verify that the function \mythetag{4.9} obeys all of the 
gluing conditions \mythetag{4.10}.\par
     Note that $h=h(t)$ in \mythetag{4.8} is a monotonic function on 
the interval $0\leqslant t\leqslant t_0$ such that $h(0)=1$ and 
$h(t_0)=0$. Therefore, we have the following lemma.
\mylemma{4.1} The function \mythetag{4.9} satisfies the inequalities 
\mythetag{4.3} for $\alpha=2$ if and only if its parameter $\varepsilon$ 
satisfies the inequalities $0\leqslant\varepsilon\leqslant 1$.
\endproclaim
     Having constructed the function $g(t)$ by means of the formula 
\mythetag{4.9}, we define the function $q=q(t)$ by applying the inverse 
conversion formula \mythetag{4.2}. It yields 
$$
\hskip -2em
q(t)=\cases 12\,t\,(1-\varepsilon\,r(t))&\text{for \ }0<t<t_0,\\
\kern 2.5em 12\,t&\text{for \ }t\geqslant t_0,
\endcases
\mytag{4.11}
$$
where $r(t)$ is a polynomial of the degree $4$. In order to simplify the 
formula for the polynomial $r(t)$ we write this polynomial as 
$$
\hskip -2em
r(t)=R(\tau)\text{, \ where \ }\tau=\frac{t_0-t}{t_0}. 
\mytag{4.12}
$$
Then the polynomial $R(\tau)$ in \mythetag{4.12} is given by the following 
formula:
$$
\hskip -2em
R(\tau)=21\,\tau^4-34\,\tau^3+16\,\tau^2-2\,\tau. 
\mytag{4.13}
$$
As we see, the polynomial \mythetag{4.13} factorizes into the product of two
polynomials:
$$
\hskip -2em
R(\tau)=(21\,\tau^3-34\,\tau^2+16\,\tau-2)\,\tau. 
\mytag{4.14}
$$
Let's denote through $R_3(\tau)$ the first multiplicand in \mythetag{4.14}. 
Then
$$
\pagebreak
\hskip -2em
R_3(\tau)=21\,\tau^3-34\,\tau^2+16\,\tau-2.
\mytag{4.15}
$$\par
\parshape 14 0pt 360pt 0pt 360pt 180pt 180pt 180pt 180pt 
180pt 180pt 180pt 180pt 180pt 180pt 180pt 180pt 180pt 180pt 
180pt 180pt 180pt 180pt 180pt 180pt 180pt 180pt 0pt 360pt 
The function \mythetag{4.15} is a cubic polynomial of the variable $\tau$. 
\vadjust{\vskip 5pt\hbox to 0pt{\kern 
10pt \includegraphics{khab03_01.eps}\hss}\vskip -5pt}
The graph of this polynomial on the segment $0\leqslant\tau\leqslant 1$ is 
shown in Fig\.~4.1. At the ending points of this segment we have 
$$
R_3(0)=-2,\qquad R_3(1)=1.\quad 
\mytag{4.16}
$$
The values \mythetag{4.16} are easily verified by means of direct
calculations. Apart from \mythetag{4.16}, there are two local extrema 
of the function $R_3(\tau)$ within this segment:
$$
\aligned
\tau_{\sssize\text{max}}=\frac{34-2\,\sqrt{37}}{63}\approx 0.34,\\
\tau_{\sssize\text{min}}=\frac{34+2\,\sqrt{37}}{63}\approx 0.73.
\endaligned
\quad
\mytag{4.17}
$$
Substituting the quantities \mythetag{4.17} into the formula 
\mythetag{4.15}, we easily find the values of the polynomial 
$R_3(\tau)$ at its local extrema $\tau_{\sssize\text{max}}$ and 
$\tau_{\sssize\text{min}}$:
$$
\hskip -2em
\aligned
&R_{\sssize\text{max}}=\frac{394+592\,\sqrt{37}}{11907}\approx
0.33,\\
&R_{\sssize\text{min}}=\frac{394-592\,\sqrt{37}}{11907}\approx
-0.26.
\endaligned
\mytag{4.18}
$$
From \mythetag{4.16} and \mythetag{4.18} we immediately derive the 
inequality
$$
\hskip -2em
R_3(\tau)\leqslant 1\text{\ \ for all \ }0\leqslant\tau\leqslant 1.
\mytag{4.19}
$$
Since $\tau\geqslant 0$ in \mythetag{4.19}, we can multiply both sides
of \mythetag{4.19} by $\tau$. This yields $R_3(\tau)\,\tau\leqslant\tau$.
Then, combining this inequality with $\tau\leqslant 1$ and taking into
account the formulas \mythetag{4.14} and \mythetag{4.15}, we obtain
$$
\hskip -2em
R(\tau)\leqslant 1\text{\ \ for all \ }0\leqslant\tau\leqslant 1.
\mytag{4.20}
$$
Applying \mythetag{4.12}, we transform \mythetag{4.20} to the inequality
$$
\hskip -2em
r(t)\leqslant 1\text{\ \ for all \ }0\leqslant t\leqslant t_0,
\mytag{4.21}
$$
while \mythetag{4.16} combined with \mythetag{4.15}, \mythetag{4.14} 
and \mythetag{4.12} yields 
$$
\xalignat 2
&\hskip -2em
r(0)=1,&&r(t_0)=0.
\mytag{4.22}
\endxalignat
$$
On the base of \mythetag{4.11}, \mythetag{4.21} and \mythetag{4.22} we 
can formulate the following lemma similar to the lemma~\mythelemma{4.1}.
\mylemma{4.2} The function \mythetag{4.11} satisfies the inequalities 
$q(t)\geqslant 0$ for all $t>0$ if and only if its parameter $\varepsilon$ 
satisfies the inequalities $0\leqslant\varepsilon\leqslant 1$.
\endproclaim
\mytheorem{4.1} For each particular value of the parameter $\varepsilon$ 
satisfying the inequa\-lities $0<\varepsilon\leqslant 1$ the function 
\mythetag{4.11} is a counterexample to Khabibullin's 
conjecture~\mytheconjecture{1.1} for \pagebreak $n=2$ and $\alpha=2$. 
\endproclaim
     The theorem~\mythetheorem{4.1} follows from the lemma~\mythelemma{4.2}
and the above considerations in \S\,3 and in \S\,4 preceding it. Indeed, 
the function \mythetag{4.11} is derived from the function \mythetag{4.9}
by means of the inverse conversion formula \mythetag{4.2}. Therefore,
substituting \mythetag{4.11} into the direct conversion formula 
\mythetag{4.1}, we get the function \mythetag{4.9}.\par
      The formula \mythetag{4.1} is related to the inequality \mythetag{1.1}
in Khabibullin's conjecture~\mytheconjecture{1.1} through the formulas
\mythetag{2.3}, \mythetag{2.2}, and \mythetag{2.1}. Therefore, substituting
the function \mythetag{4.11} into the left hand side of the inequality 
\mythetag{1.1} for $n=2$, we get
$$
\hskip -2em
\int\limits^{\,1}_0\left(\,\,\int\limits^{\,1}_x(1-y)^{n-1}
\,\frac{dy}{y}\right)q(tx)\,dx=\frac{g(t)}{t},
\mytag{4.23}
$$      
where $g(t)$ is given by the formula \mythetag{4.9}. Now let's recall that
$\alpha=2$ and apply the lemma~\mythelemma{4.1} to \mythetag{4.23}. This
yields the following result.
\mylemma{4.3}For each particular value of the parameter $\varepsilon$ 
such that\/ $0\leqslant\varepsilon\leqslant 1$ the function \mythetag{4.11} 
satisfies the inequality \mythetag{1.1} for $n=2$ and $\alpha=2$. 
\endproclaim
     The next step is to calculate the integral in the left hand side of the 
inequality \mythetag{1.2} for the function \mythetag{4.11}. For this purpose
we use the formula 
$$
\hskip -2em
\int\limits^{+\infty}_0\!\!q(t)\,\ln\Bigl(1+\frac{1}
{t^{\,2\kern 0.2pt\alpha}}\Bigr)\,dt\,=\int\limits^{+\infty}_0
\!\!\varPhi_{n-1}(\alpha,t)\,g(t)\,dt.
\mytag{4.24}
$$
The formula \mythetag{4.24} is derived from \mythetag{4.1} with the use 
of the formula \mythetag{2.7}. Let's recall that in our case $\alpha=2$ 
and $g(t)$ is given by the formula \mythetag{4.9}. Then the formula 
\mythetag{4.24} is written in the following form:
$$
\int\limits^{+\infty}_0\!\!q(t)\,\ln\Bigl(1+\frac{1}
{t^{\,2\kern 0.2pt\alpha}}\Bigr)\,dt\ =\!\int\limits^{+\infty}_0
\!\!\varPhi_{n-1}(\alpha,t)\ t^\alpha\,dt
\,-\,\varepsilon\!\!\int\limits^{\ t_0}_0
\!\!\varPhi_{n-1}(\alpha,t)\,t^2\,h(t)\,dt.
\qquad
\mytag{4.25}
$$
Now we apply the formula \mythetag{3.13} to \mythetag{4.25} and write this 
formula as
$$
\hskip -2em
\int\limits^{+\infty}_0\!\!q(t)\,\ln\Bigl(1+\frac{1}
{t^{\,2\kern 0.2pt\alpha}}\Bigr)\,dt=\pi\,\alpha\prod^{n-1}_{k=1}
\Bigl(1+\frac{\alpha}{k}\Bigr)
+\,\delta I.
\mytag{4.26}
$$
The term $\delta I$ in the right hand side of \mythetag{4.26} is given by 
the formula 
$$
\hskip -2em
\delta I=-\,\varepsilon\!\int\limits^{\ t_0}_0\!\varPhi_{n-1}(\alpha,t)
\ t^2\,h(t)\,dt.
\mytag{4.27}
$$\par
    In our case the parameter $\varepsilon$ is positive due to the condition
$0<\varepsilon\leqslant 1$ in the theorem~\mythetheorem{4.1}. The function
$h(t)$ is given by the formula \mythetag{4.8}. It is positive for
$0<t<t_0$. As for the function $\varPhi_{n-1}(\alpha,t)$, since $n=2$ and
$\alpha=2$, in our case $\varPhi_{n-1}(\alpha,t)$ is given by the formula
\mythetag{4.6}. It is negative for $0<t<t_0$. As a result we derive
the following inequality for the integral \mythetag{4.27}:
$$
\hskip -2em
\delta I>0
\mytag{4.28}
$$
The inequality \mythetag{4.28} completes the proof of the 
theorem~\mythetheorem{4.1}.\par 
    The integral \mythetag{4.27} can be computed
numerically. The product in the right hand side of \mythetag{4.27} can
also be computed numerically. As a result we get
$$
\hskip -2em
\int\limits^{+\infty}_0\!\!q(t)\,\ln\Bigl(1+\frac{1}
{t^{\,2\kern 0.2pt\alpha}}\Bigr)\,dt\approx 
18.84955592+0.01299443\,\varepsilon. 
\mytag{4.29}
$$
On the other hand there is the estimate \mythetag{3.15} for the integral
\mythetag{4.29}. The constant $C(n,\alpha)$ in \mythetag{3.15} is given 
by the integral \mythetag{3.10}. It can be calculated numerically. For 
$n=2$ and $\alpha=2$ we have $C(n,\alpha)\approx 19.65507202$. Then 
\mythetag{3.15} yields
$$
\hskip -2em
\int\limits^{+\infty}_0\!\!q(t)\,\ln\Bigl(1+\frac{1}
{t^{\,2\kern 0.2pt\alpha}}\Bigr)\,dt\leqslant 
19.65507203.
\mytag{4.30}
$$
Note that $0<\varepsilon\leqslant 1$ in the theorem~\mythetheorem{4.1}. 
As we see, even for $\varepsilon=1$ there is a substantial gap between 
the integral \mythetag{4.29} and the estimate \mythetag{4.30} for it.  
\head
5. Conclusions. 
\endhead
     The theorem~\mythetheorem{4.1} shows that Khabibullin's conjecture 
is not valid for $\alpha>1/2$ in its present form~\mytheconjecture{1.1}. 
However, due to the theorem~\mythetheorem{3.1} and the inequalities 
\mythetag{3.14} it can be replaced by the following theorem which is 
certainly valid for all $\alpha>0$. 
\mytheorem{5.1} For each positive integer $n>0$ and for each positive 
$\alpha>0$ there is a positive constant $C[\Kh](n,\alpha)$ that yields 
the best (non-improvable) estimate 
$$
\hskip -2em
\int\limits^{+\infty}_0\!\!q(t)\,\ln\Bigl(1+\frac{1}
{t^{\,2\kern 0.2pt\alpha}}\Bigr)\,dt\leqslant C[\Kh](n,\alpha)
\mytag{5.1}
$$
in the class of all non-negative continuous functions 
$q(t)\geqslant 0$ on the half-line $t>0$ satisfying the integral 
inequality 
$$
\hskip -2em
\int\limits^{\,1}_0\left(\,\,\int\limits^{\,1}_x(1-y)^{n-1}
\,\frac{dy}{y}\right)q(tx)\,dx\leqslant t^{\alpha-1}
\text{\ \ for all \ }t>0.
\mytag{5.2}
$$
\endproclaim
     In \mythetag{5.2} we omit the inessential value $t=0$ as compared 
to the initial inequality \mythetag{1.1} in Khabibullin's 
conjecture~\mytheconjecture{1.1}.\par
     Note that the constants $C[\Kh](n,\alpha)$ in \mythetag{5.1} are 
finite. They satisfy the inequality $C[\Kh](n,\alpha)\leqslant 
C(n,\alpha)<\infty$, where the constants $C(n,\alpha)$ are given by 
the formula \mythetag{3.10}. I think the constants $C[\Kh](n,\alpha)$ 
in \mythetag{5.1} should be called the Khabibullin constants, appreciating 
the efforts of Prof\. B\.~N\.~Khabibullin in formulating the 
conjecture~\mytheconjecture{1.1} and advertising this conjecture for
many years.\par
     Note that the theorem~\mythetheorem{5.1} claims the existence of 
the constants $C[\Kh](n,\alpha)$, but it gives neither a formula nor 
an algorithm for calculating these constants. Nowadays the problem of 
finding the exact values of the  Khabibullin constants $C[\Kh](n,\alpha)$ 
in \mythetag{5.1} is yet an unsolved problem. 

\enddocument
\end